\documentclass[12pt]{article}

\usepackage[cp866]{inputenc}       % Кодировк  входного документ ;
\usepackage{amsfonts}

\usepackage[T2A]{fontenc}           % Кодировк  для шрифтов LH

\usepackage[english]{babel} % Включение русифик ции, русских и
\usepackage{epsfig}
\usepackage{graphicx}
\catcode`\@=11
\@addtoreset{equation}{section}\catcode`\@=12

\begin{document}

%\documentstyle[11pt,epsfig]{article}
%\textwidth 17.0 cm
%\textheight 9.5in
%\topmargin -1in
%\oddsidemargin -0.5in
%\evensidemargin -0.5in
%\begin{document}

\newtheorem{lm}{Lemma}
\newtheorem{theorem}{Theorem}
\newtheorem{df}{Definition}
\newtheorem{prop}{Proposition}

\begin{center}
 {\large\bf Homoclinic tangencies to resonant saddles and discrete Lorenz attractors}
\vspace{12pt}

 {\bf S.V.Gonchenko}$^1$, {\bf I.I.Ovsyannikov}$^{12}$
 \label{Author}
\vspace{6pt}

 {\small
  $^1$ Nizhny Novgorod State University, Russia; \\
  $^2$ Universit\"at Bremen, Germany \\

  E-mail: sergey.gonchenko@mail.ru; ivan.i.ovsyannikov@gmail.com \\
 }

\end{center}

{\bf Abstract.}
We study bifurcations of periodic orbits in three parameter general unfoldings of
certain types quadratic homoclinic tangencies to saddle fixed points. We apply the rescaling
technique to first return (Poincar\'e) maps and show that the rescaled maps can be brought to a map
asymptotically close to the 3D Henon map $\bar x=y,\bar y=z,\bar z = M_1 + M_2 y + B x - z^2$
which, as known \cite{GOST05}, exhibits wild hyperbolic Lorenz-like attractors in some open domains
of the parameters. Based on this, we prove the existence of infinite cascades of Lorenz-like
attractors\footnote{This work was supported by grant 14-41-00044 of the RSF and
grant of RFBR 13-01-00589.}.


\\~\\

\section{Introduction}
\begin{figure}[ht]
\label{fig:Lorenz}
%\begin{tabular}{cc}
%\hspace*{-8mm}\psfig{file=figGST1a4.ps,width=70mm,angle=-90} &
%\hspace*{-23mm}\psfig{file=figGST1b.ps,width=70mm,angle=-90}
%\end{tabular}
\centerline{\epsfig{file=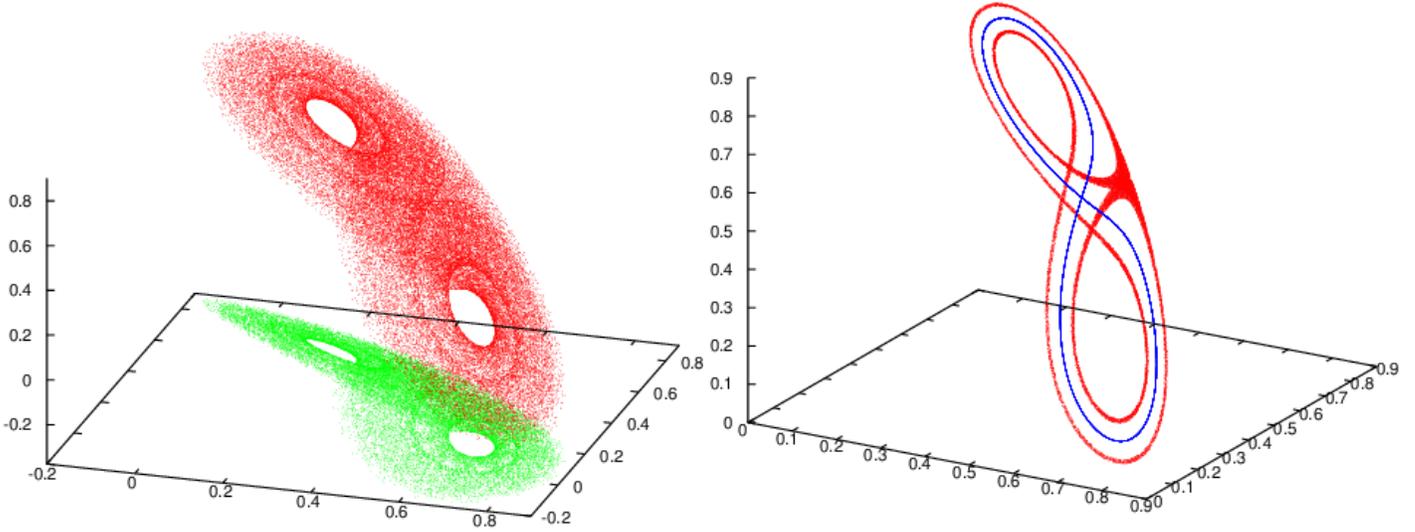, width=20cm}}
\caption{
{\footnotesize Plots of attractors of map (\ref{3Hen1}) observed numerically in \cite{GOST05} for
$M_1=0,B=0.7$ and $M_2=0.85$ (left) or $M_2=0.815$ (right). In
the left panel, the projection on the $(x,y)$-plane is also displayed.
In the right panel, a  ''figure-eight'' saddle closed invariant
curve inside the lacuna is shown. Note the similarity to the Lorenz
attractors of the Shimizu-Morioka system, see \cite{ShA86,ShA89,ShA93}.}} 
\end{figure}

In \cite{GOST05} it was discovered that the three-dimensional Henon map
\begin{equation}
\bar x = y, \; \bar y = z, \; \bar z = M_1 + M_2 y + Bx - z^2,
\label{3Hen1}
\end{equation}
where $(M_1,M_2,B)$ are parameters ($B$ is the Jacobian of map),
can possess strange attractors that seem
very similar to the Lorenz attractors, see fig.~\ref{fig:Lorenz}. Later it was shown that such {\em discrete Lorenz attractors} can arise as result of simple, universal and natural bifurcation scenarios realizing in one-parameter families of three-dimensional maps \cite{GGS12,GGKT14}. This means, in fact, that the discrete Lorenz attractors can be met widely in applications. For instance, 
%As a confirmation of this, we refer the reader to 
in \cite{GG14,GG15} such attractors were found in nonholonomic models of rattleback (called also as a Celtic stone). See also \cite{GGOT13,GG15} where various types of strange homoclinic attractors, including discrete Lorenz ones, were investigated.

The similarity between the discrete and classical Lorenz attractors appears to be not accidental and it can be explained by various reasons. Thus, it is well known that the classical Lorenz attractor can be born as a result of local bifurcations of an equilibrium state with three zero eigenvalues when a flow possesses a 
%additional conditions including those related to the existence of some
(Lorenzian) symmetry \cite{SST93}. Analogously for maps,  discrete Lorenz attractors can arise under bifurcations of fixed points with multipliers $(-1, -1, +1)$, in this case the required local symmetry exists automatically due to negative multipliers. As it was shown in \cite{GOST05, GGOT13}, the second iteration of the map near this point can be embedded into a flow up to asymptotically small periodic non-autonomous terms. The corresponding flow normal form of such bifurcations coincides with the well-known Shimizu-Morioka model, which, in turn, exhibits the Lorenz attractor for certain parameter values \cite{ASh93,TT11}. Thus, we can consider the attractor in the map as the one of the Poincar\'e map (period map) of a periodically perturbed system with the Lorenz attractor. On the other hand, as it was shown in  paper \cite{TS08} by Turaev and Shilnikov, such discrete attractor is genuine in the sense that every its orbit has positive maximal Lyapunov exponent\footnote{Moreover, it is a wild pseudohyperbolic attractor \cite{TS08}, since it allows homoclinic tangencies (i.e. contains Newhouse wild hyperbolic sets \cite{N79}) and has an adsorbing domain inside which the differential $DT$, for the map itself and all close maps can be decomposed into a direct sum of transverse invariant subspaces ${\cal W}^{ss}$ and ${\cal W}^{eu}$ where $DT_{{\cal W}^{ss}}$ is strongly contacting and $DT_{{\cal W}^{eu}}$ expands exponentially volume, i.e.  $\|DT^k_{{\cal W}^{ss}}\|< L \sigma^k$ and $|\det{DT^k_{{\cal W}^{eu}}}|> L \nu^k$ for some constants $L > 0, 0 < \sigma < 1 < \nu$ and all positive $k$.}.

In the present paper we study bifurcations of three-dimensional diffeomorphisms with homoclinic tangencies, leading to the birth
of descrete Lorenz attractors. Problems of this kind were previously analyzed in \cite{GMO06, GST09, GO10, GO13, GOT14}.

In \cite{GST09, GO10, GO13} the birth of discrete Lorenz attractors from nontransversal heteroclinic cycles of three-dimensional diffeomorphisms was studied.
Such a cycle contains two fixed points $O_1$ and $O_2$ of type (2,1), i.e. with $\dim W^s(O_i)=2, \dim W^u(O_i)=1$, and one pair of stable and unstable manifolds intersect transversely and another pair has a quadratic tangency. %For more definiteness, let $W^u(O_1)$ and $W^s(O_2)$ intersect transversely at the points of a heteroclinic orbit $\Gamma_{12}$ and $W^u(O_2)$ and $W^s(O_1)$ have a quadratic tangency at the points of a heteroclinic orbit $\Gamma_{21}$. 
It was assumed that at least one of the points $O_1$ and $O_2$ is a saddle-focus, see Fig. Moreover, in all cases the additional condition that
%HoweverGO13},  e and such that at least one of them is a saddle-focus and
the Jacobians
of the map in points $O_1$ and $O_2$ are greater and less than one respectively was imposed (the so-called case of contracting-expanding maps). The birth of discrete Lorenz attractors was proved for three-parameter general unfoldings. \\

{\bf Remark.} Naturally, three parameters are needed to allow generically the existence of triply degenerate fixed points in the corresponding first return maps. In such families the first return map
can be rescaled to the form asymptotically close to map (\ref{3Hen1}).
Thus, using the results of \cite{GOST05} (see also \cite{GGOT13} for more generic statement),
we deduce the birth of discrete Lorenz attractors in close systems. \\

In the case of homoclinic tangencies to the saddle fixed point $O$ of a three-dimensional diffeomorphism $T$ the birth of Lorenz attractors was proved in the cases when:

1) \cite{GMO06}, the point $O$ is a saddle-focus with the unit Jacobian (saddle-focus of conservative type). % , see fig.~\ref{fig:glob}a.

2) \cite{GOT14}, the fixed point is a saddle with the unit Jacobian and the quadratic tangency is non-simple\footnote{The definition of simple homoclinic
tangency can be found in \cite{GST93c}. In particular, it assumes that the so-called extended unstable invariant manifold intersects transversely
the leaf of the strong stable foliation in the point of tangency. The main cases of non-simple homoclinic tangencies in three-dimensional
diffeomorphisms were considered in \cite{Tat01}, see also condition {\bf D} in \S 2 of the present paper.}. %, see fig.~\ref{fig:glob}b.

We note that the condition on Jacobians in all these cases is necessary for the existence of a non-trivial (three-dimensional) dynamics
in the neighborhood of the homoclinic orbit \cite{T96}. Otherwise, if, for example, one has $J < 1$, all three-dimensional volumes will
be contracted under the iterations of map $T$ (near point $O$) and, hence, the dynamics of first return maps  $T^k$ for
large $k$ will be effectively two-dimensional, or even one-dimensional. Recall that, by definition \cite{T96}, the {\em effective dimension} $d_e$ of a bifurcation problem equals $n$, if periodic orbits with $n$ multipliers equal $\pm 1$ can appear at bifurcations but no orbits exist with more than $n$ unit multipliers.

Formally, if to consider three-dimensional diffeomorphisms with homoclinic tangencies to a hyperbolic saddle fixed point with $|J|=1$, then $d_e$ can be equal to $3$ only in the following cases: (i) the point is a saddle-focus; (ii) the point is the saddle (all multipliers are real) and the tangency is not simple, and (iii) the point is a resonant saddle, i.e. it has two multipliers equal in the absolute value. Otherwise, the effective dimension is less than three since the direction of strong contraction is present \cite{GST93c,GST08} for all nearby systems.

Cases (i) and (ii) were considered in \cite{GMO06} and \cite{GOT14} respectively.
In this paper we consider the new case (iii) when the saddle is resonant. Note that if the resonance $\lambda_1=\lambda_2$ takes place, we may perturb the map in such a way that the resulting map will have a saddle-focus fixed point with $|J|=1$ and, hence, we can apply results of \cite{GMO06} to prove the birth of discrete Lorenz attractors.
It not the case for the resonance $\lambda_1= - \lambda_2$ which is of independent interest.

We consider the case when a fixed
point $O$ has multipliers $\lambda, -\lambda, \gamma$ such that $0 < \lambda < 1$, $|\gamma| > 1$ and $|\lambda^2 \gamma| = 1$.
This means that $O$ is a resonant saddle point of conservative type. Obviously, the bifurcation codimension of this problem is at least three and, as we will show, $d_e=3$ in this case.

We show that in the three-parametric families $f_\mu$, $\mu = (\mu_1, \mu_2, \mu_3)$ unfolding generally this type of a homoclinic tangency,
in the parameter space there exist domains $\triangle_k \to \{\mu = 0\}$ as $k \to \infty$ such that for $\mu \in \triangle_k$
the first return map $T_k$ possesses the discrete Lorenz attractor. Recall that the map $T_k$ is constructed by the iterations of  map $f_\mu$, i.e. $T_k=f_\mu^k$, but the domain of its definition is a small box $\sigma_0^k$ near some homoclinic point.

The paper consists of two paragraphs. In \S 2 we formulate our main result -- Theorem~\ref{th3-2} and construct the first return map
of some small neighborhood of the homoclinic orbit. In \S 3 we prove Theorem~\ref{th3-2}.

\iffalse

In the paper we consider one more new case of a homoclinic tangency to a fixed point of the saddle type. We assume that
point $O$ has multipliers $\lambda, -\lambda, \gamma$ such that $0 < \lambda < 1$, $\gamma > 1$ and $|\lambda^2 \gamma| = 1$.
This means that $O$ is a resonant saddle point of conservative type. Obviously, the codimension of this problem is three.

We show that in the three-parametric families $f_\mu$, $\mu = (\mu_1, \mu_2, \mu_3)$ unfolding generally this type of a homoclinic tangency,
in the parameter space there exist domains $\triangle_k \to \{\mu = 0\}$ as $k \to \infty$ such that for $\mu \in \triangle_k$
the first return map $T^k$ possesses the discrete Lorenz attractor.

The paper consists of two paragraphs. In \S 2 we formulate our main result -- Theorem~\ref{th3-2} and construct the first return map
of some small neighborhood of the homoclinic orbit. In \S 3 we prove Theorem~\ref{th3-2}.

\fi

\section{Statement of the problem and formulation of main results.}

We study bifurcations of three-dimensional diffeomorphisms of a special type
(codimension two) quadratic homoclinic tangency to a saddle fixed point with the unit Jacobian.
Namely, we assume that the initial diffeomorphism $f_0 \in C^r, r\geq 5,$ satisfies the following
conditions:

\begin{enumerate}
\item[{\bf A)}] $f$ has a saddle fixed point $O$ with real multipliers $\lambda_1, \lambda_2, \gamma$ such
that  $0 < |\lambda_{1, 2}| < 1 < |\gamma|$ and
$$
J_0 \equiv |\lambda_1 \lambda_2 \gamma| = 1.
$$
\item[{\bf B)}] The stable $W^s(O)$ and unstable $W^u(O)$ invariant manifolds of $O$ have a quadratic
tangency at the points of some homoclinic orbit $\Gamma_0$.
\item[{\bf C)}]  The saddle $O$ is resonant in the sense that $\lambda_1 = -\lambda_2 = \lambda > 0$.
\end{enumerate}

Condition {\bf A} means that the point $O$ is a saddle of conservative type and $\dim \; W^s(O) = 2$
and $\dim \; W^u(O) = 1$. Condition {\bf C} is an additional degeneracy of the saddle fixed point.
We will consider smooth parameter families $f_\varepsilon$ of
diffeomorphisms (general unfoldings of conditions
{\bf A}--{\bf C}), such that $f_0$ belongs to it for $\varepsilon = 0$.

Let $U\equiv U(O\cup\Gamma_0)$ be a sufficiently small fixed
neighbourhood of $\Gamma_0$  that is a union
of a neighbourhood $U_0$ of  $O$ and a number of
neighbourhoods of those points of $\Gamma_0$ which lie outside $U_0$.
Denote by  $T_{0}$ the restriction of the diffeomorphism $f_\varepsilon$ onto $U_0$.
We call $T_{0}$ {\it a local map.} By a linear transformation of
coordinates in $U_0$, map $T_0$ can be written as
%in the form
$$
(\bar x_1,\bar x_2,\bar y) = (\lambda_1 x_1,\lambda_2 x_2,\gamma y)\;
+\; h.o.t.
$$
The origin $O=(0,0,0)$ is a fixed point of $T_{0}$, the stable
manifold $W^{s}(O)$ is tangent at $O$ to the $(x_1,x_2)$-plane and
the unstable manifold $W^{u}(O)$ is tangent at $O$ to the
$y$-axis. The intersection points of $\Gamma_{0}$ with $U_0$ belong to
the set $W^{s}\cap W^{u}$ and accumulate to $O$ at both forward and backward iterations. Thus, infinitely
many points of $\Gamma_0$ lie on $W^{s}_{loc}$ and $W^{u}_{loc}$.
Let $M^{+} \in  W^{s}_{loc}$   and $M^{-} \in W^{u}_{loc}$ be two
such points and let $M^+ = f_0^{n_0}(M^-)$ for some positive
integer $n_0$. Let $\Pi^+\subset U_0$  and $\Pi^-\subset U_0$ be
small neighbourhoods of points $M^{+}$  and $M^{-}$ respectively.
The map $T_{1}\equiv f_\varepsilon^{n_0}: \Pi^- \rightarrow  \Pi^+$ is called
{\it a global map.}
%Denote $L_u = W^u_{loc}\cap\Pi^-$ and $L_s =
%W^s_{loc}\cap\Pi^+$.

%It follows from condition B) that $T_1(W^u_{loc})$ and $W^s_{loc}$ have a
%quadratic tangency at $M^+$. Then, by the definition, we can
%choose local $C^2$-coordinates $(\xi_1,\xi_2,\eta)$ near $M^+$
%such that the equations of  $W^s_{loc}$ and $T_1(W^u_{loc})$ can be written in
%the following parametric form
%\begin{equation}
%\begin{array}{cl}
%W^s_{loc}\;& :\;\;\;
%\eta=0  \;, \\
%%
%T_1(W^u_{loc})\;&:\;\;\xi_1=\alpha_1 t + O(t^2)\;,\; \xi_2=\alpha_2 t +
%O(t^2)\;,\; \eta=\beta t^2 + o(t^2)\;,
%\end{array}
%\label{eqsu}
%\end{equation}
%with parameter $t$. Then the tangency is quadratic if
%\begin{equation}
%\begin{array}{l}
%\alpha_1^2+\alpha_2^2 \neq 0\;\;{\rm and}\;\; \beta \neq 0\;.
%\end{array}
%\label{quta}
%\end{equation}

%\subsection{The main normal form of the local map $T_0$.}\label{mapt01}

From \cite{book, GS90, GS92, GST07} it is known that there exists a $C^{r}$-change of coordinates (which is
$C^{r-2}$-smooth in the parameters) bringing $T_0$ to the so-called {\em main normal form}:
%
%
%We will use in $U_0$ such local coordinates in which the $T_0$ is written in the so-called {\em main
%normal form} \cite{book, GS90, GS92, GST07}.
%There exists a $C^{r}$-change of coordinates (which is
%$C^{r-2}$-smooth in the parameters) bringing $T_0$ to the
%following form
\begin{equation}
\begin{array}{l}
\bar x_1 \; = \; \lambda_1(\varepsilon) x_1 +
O(\|x\|^2|y|) \;, \\
\bar x_2 \; = \; \lambda_2(\varepsilon) x_2 +  O(\|x\|^2|y|) \;,\\
\bar y \; = \; \gamma(\varepsilon) y + O(\|x\||y|^2) \;.\\
\end{array}
\label{t0norm}
\end{equation}
The main peculiarity of this form is that in coordinates (\ref{t0norm}) the stable and unstable manifolds of the saddle
fixed point are locally straightened, their equations are $W^s:  \{y = 0\}$, $W^u:  \{x_1 = 0, x_2 = 0\}$. The main normal form also
allows to obtain a quite simple representation of the iterations of $T_0$. The latter can be formulated as the following lemma:

\begin{lm} {\rm \cite{GST08}}
For any positive integer $k$ and for any sufficiently small $\varepsilon$ the map
$T_0^k(\varepsilon)\;:\; (x_0,y_0)\;\to\; (x_k,y_k)$ can be written in the following cross-form
\begin{equation}
\begin{array}{l}
x_{k1} \;-\; \lambda_1^{k}(\varepsilon) x_{01} \;=\;
\hat\lambda^{k}\xi_{k1}(x_0,y_k,\varepsilon) \;,\; \\
x_{k2} \;-\; \lambda_2^{k}(\varepsilon) x_{02} \;=\;
\hat\lambda^{k}\xi_{k2}(x_0,y_k,\varepsilon) \;,\; \\
y_0 \;-\; \gamma(\varepsilon)^{-k} y_k \;=\;
\hat\gamma^{-k}\eta_k(x_0,y_k,\varepsilon) \;\;,
\end{array}
\label{T0kk}
\end{equation}
where $\hat\lambda$ and $\hat\gamma$ are some constants such that $\hat\lambda = \lambda +
\delta, \; \hat\gamma = \lambda|\gamma^{-1}| -\delta$ for some small $\delta>0$  and functions $\xi_{k}$ and
$\eta_{k}$ are uniformly bounded along with all derivatives up to order $(r-2)$. \label{lem2}
\end{lm}

%\subsection{Properties of the global map $T_1$.}\label{propt1}

Next we construct the most appropriate form for the global map $T_1$ for all small $\varepsilon$.
Let the chosen homoclinic points have coordinates $M^+ = M^+(x_1^+,x_2^+,0)\in W^s_{loc}$ and $M^- = M^-(0,0,y^-)\in W^u_{loc}$,
where $(x_1^+)^2 + (x_2^+)^2\neq 0$ and $y^- > 0$.
At $\varepsilon = 0$ we have that $T_1M^- = M^+$ and $T_1(W_{loc}^u)$ and $W_{loc}^s$ are
tangent quadratically at the point $M^+$. Thus, the global map $T_1$ at $\varepsilon=0$ can
be written as the Taylor expansion near the point $(x_1 = 0, x_2 = 0, y = y^-$):
\begin{equation}
\label{t1-10}
\begin{array}{l}
\bar x_1 - x_1^+ \; = \; a_{11}x_1 + a_{12}x_2 + b_1(y-y^-) + O(\|x\|^2) + O(\|x\| |y - y^-|) + O((y - y^-)^2)
\\
\bar x_2 - x_2^+ \; = \; a_{21}x_1 + a_{22}x_2 +
b_2(y - y^-) + O(\|x\|^2) + O(\|x\| |y - y^-|) + O((y - y^-)^2) \\
\bar y \; = \; c_1x_1 + c_2x_2 +
d(y - y^-)^2 + O(\|x\|^2) + O(\|x\| |y - y^-|) + O(|y - y^-|^3)
%
%\\
\end{array}
\end{equation}
The equation of curve $T_1(W^u_{loc})$  at  $\varepsilon = 0$ looks as follows
(we put $x_1=x_2=0$ in (\ref{t1-10})):
\begin{equation}
\label{t1-11}
\begin{array}{l}
\bar x_1 - x_1^+ \; = \; b_1(y - y^-) + O((y - y^-)^2)\\
\bar x_2 - x_2^+ \; = \; b_2(y - y^-) + O((y - y^-)^2)\\
\bar y \; = \; d(y - y^-)^2 + O(|y - y^-|^3) \\
\end{array}
\end{equation}
This is a parametric equation (with parameter $(y - y^-)$) of the
curve $T_1(W^u_{loc})$ in a neighbourhood of $M^+$. The equation
of $W^s_{loc}$ is $y=0$. Since the initial homoclinic tangency is
quadratic, it follows that $ d \neq 0,\; b_1^2 + b_2^2 \neq 0 $. Moreover, map $T_1(0)$ is a
diffeomorphism, therefore
\begin{equation}
\label{t1-3}
J_1 = {\rm det}\; \left(
\begin{array}{rcl}
a_{11}\; & a_{12}\; &\;b_1 \\
a_{21}\; & a_{22}\; &\;b_2 \\
c_1\; & c_2 &\; 0 \\
\end{array}
\right)
\neq 0
\end{equation}
and, hence, $c_1^2 + c_2^2\neq 0$.

At small $\varepsilon$ the global map
$T_1(\varepsilon)$ can be written in the following form
(the Taylor expansion near the point $(x_1,x_2,y) = (0,0,y^-(\varepsilon))$
\begin{equation}
\label{t1-1}
\begin{array}{l}
\bar x_1 - x_1^+(\varepsilon) \; = \; a_{11}x_1 + a_{12}x_2 + b_1(y-y^-(\varepsilon))  + \ldots
\\
\bar x_2 - x_2^+(\varepsilon) \; = \; a_{21}x_1 + a_{22}x_2 +
b_2(y-y^-(\varepsilon)) + \ldots \\
\bar y \; = \; y^+(\varepsilon) + c_1x_1 + c_2x_2 + d(y-y^-(\varepsilon))^2 + \ldots
\\
\end{array}
\end{equation}
where $x_1^+(0)=x_1^+,x_2^+(0)=x_2^+,y^-(0)= y^-,y^+(0)= 0$;
%and
all coefficients $a_{11}, \ldots, d$ depend (smoothly) on $\varepsilon$; and we shift $y^-$ into
$y^-(\varepsilon)$ in order to nullify the linear in $y$ terms from the right side of the third
equation.

We assume that the following general condition  holds
\begin{equation}
 b_1 c_1 b_2 c_2 \neq 0.
\label{b1c1n0}
\end{equation}
Its meaning is as follows: if at least one of these four coefficients is zero, then in any neighborhood of $f_0$ there exist maps having a nontransversal homoclinic orbit
close to $\Gamma_0$ with a non-simple quadratic tangency (see \cite{Tat01, GGT07, GOT14} for details); it composes an additional degeneracy which
we do not consider here.

Conditions {\bf A}, {\bf B}, {\bf C} together with (\ref{b1c1n0}) define a codimension $3$ bifurcation surfaces of diffeomorphisms with a quadratic homoclinic tangency. Hence, as a general unfolding we should consider a three-parameter families where the parameters $\mu_1,\mu_2$ and $\mu_3$ control the degeneracies imposed due to conditions {\bf B},
{\bf A} and {\bf C}, respectively.

%We  consider three-parameter families $f_{\mu_1,\mu_2,\mu_3}$ of diffeomorphisms close to $f$.
Naturally, the splitting distance of manifolds $W^s(O)$ and $W^u(O)$ with respect to
the point $M^+$ is considered as the first governing parameter $\mu_1$. It is seen from (\ref{t1-1}) that
\begin{equation}
\;\;\mu_1 \equiv y^+(\varepsilon)\;.
\label{mu1}
\end{equation}

The second parameter should control the Jacobian $J = \lambda_1\lambda_2\gamma$ of $f_\mu$ at
saddle $O_\mu$. Therefore, we define
\begin{equation}
\;\mu_2 = 1 - |\lambda_1\lambda_2\gamma|
\label{mu2}
\end{equation}

As the third parameter $\mu_3$ we consider the value that controls the difference between
$|\lambda_1|$ and $|\lambda_2|$, namely:
\begin{equation}
\;\mu_3 = \frac{|\lambda_1(\varepsilon)|}{|\lambda_2(\varepsilon)|} - 1
\label{mu3}
\end{equation}

Thus, the family $f_{\mu_1,\mu_2,\mu_3}$ constructed above can be considered as a general unfolding of the
corresponding homoclinic tangency to a resonant saddle, satisfying conditions {\bf A}, {\bf B} and {\bf C}.\\

Now we are able to construct the first return maps $T_k$ using formulae (\ref{T0kk}) and
(\ref{t1-1}). As a result we will obtain a formula for $T_k$ in the initial (small) variables
$(x_1, x_2, y) \in U_0$ and parameters $\mu_1$, $\mu_2$ and $\mu_3$. Next, we rescale the initial variables and
parameters
$$
(x_1, x_2, y) \mapsto (X_1, X_2, Y)\;,\; (\mu_1, \mu_2, \mu_3) \mapsto (M_1, M_2, M_3)\;,
$$
with asymptotically small (as $k \to \infty$) factors, in such a way that in
the rescaled variables and parameters map $T_k$ is rewritten as some three-dimensional quadratic map which contains
asymptotically small (as $k \to \infty$) terms. Moreover, new coordinates $(X_1, X_2 ,Y)$ and
parameters $(M_1, M_2, M_3)$ can take arbitrary finite values at large $k$ (i.e. covering all values in the limit
$k \to \infty$).

Our main result is the following theorem.

\begin{theorem}
Let $f_{\mu_1, \mu_2, \mu_3}$ be the family under consideration. Then, in the
$(\mu_1, \mu_2, \mu_3)$-parameter space, there exist infinitely many regions $\Delta_k$
ac\-cu\-mu\-la\-ting at the origin as $k \to \infty$ such that the map $T_k$  in appropriate rescaled
coordinates and parameters is asymptotically $C^{r - 1}$-close to the following limit map

\begin{equation}
\label{tk-itog}
\begin{array}{l}
\displaystyle \bar X_1 \; = \; Y,\; \;
\bar X_2  \; = X_1 ,\; \;
\bar Y
 = M_1 + M_2 X_1 + B X_2 - Y^2 \;,
\end{array}
\end{equation}
where
\begin{equation}
\label{M1-itog}
\begin{array}{l}
M_1\; = \; -d\gamma^{2k}[\mu_1 + \lambda_1^{k}c_1x_1^+ + \lambda_2^{k}c_2x_2^+ + o(\lambda^k)]
\end{array}
\end{equation}
and
\begin{equation}
\label{M2-itog}
\begin{array}{l}
\displaystyle M_2 = \left( b_1c_1 + b_2c_{2} \frac{\lambda_2^{k}}{\lambda_1^{k}}
\right)\lambda_1^k\gamma^k(1+\dots),\;\; B = J_1 (\lambda_1 \lambda_2 \gamma)^k(1 + \dots)
\end{array}
\end{equation}
\label{th3-2}
\end{theorem}

%%%%%%%%%%%%%%%%%%%%%%%%%%%%%%%%%%%%%%%%%%%%%%%%%%%%%%%%%%%%%%%%%%%%%%
%% Remark on 3D Henon maps and discrete Lorenz attractors
%%%%%%%%%%%%%%%%%%%%%%%%%%%%%%%%%%%%%%%%%%%%%%%%%%%%%%%%%%%%%%%%%%%%%%

\section{Proof of Theorem~\ref{th3-2}.}\label{sec:th3proof}

%Recall that $|\lambda_1|>|\lambda_2|$, $|\lambda_1\gamma|>1$ and $|\lambda_2\gamma|>1$ under the
%assumption A;  also we take $\hat\gamma = |\lambda_1^{-1}\gamma|$ due to lemma \ref{lem3} . Note
%that map $T_k$ is rescaled differently in the cases I and II. However, there is a preparation part
%of the proof that is conducted in the same way for both the cases.

%\subsection{Preparation form of map $T_k$ for rescaling.}\label{prform}

Using (\ref{t1-1}) and (\ref{T0kk}) one can write the map $T_k =
T_1T_0^k$ for sufficiently large $k$  and small $\varepsilon$ in
the form
\begin{equation}
\label{tk-1}
\begin{array}{l}
\bar x_1 - x_1^+ \; = \; a_{11}(\lambda_1^{k} x_{1}
+ \hat\lambda^{k}\xi_{k1}(x, y, \varepsilon)) +
a_{12}(\lambda_2^{k} x_{2} +
\hat\lambda^{k} \xi_{k2}(x, y, \varepsilon)) + \\
\qquad  + b_1(y - y^-) + O(|y - y^-|^2 + \lambda^{k} \|x\| |y - y^-|) +
\lambda^{2k}\|x\|^2)\;,\\ \\
%\\
\bar x_2 - x_2^+ \; = \;
a_{12}(\lambda_1^{k} x_{1}
+ \hat\lambda^{k}\xi_{k1}(x, y, \varepsilon)) +
a_{22}(\lambda_2^{k} x_{2} +
\hat\lambda^{k}\xi_{k2}(x, y, \varepsilon)) + \\
\qquad + b_2(y - y^-) + O(|y - y^-|^2 + \lambda^{k} \|x\| |y - y^-|) +
\lambda^{2k} \|x\|^2)\;,\\ \\
\gamma^{-k} \bar y  -
\hat\gamma^{-k}\eta_k(\bar x, \bar y, \varepsilon) \; = \;
\mu_1 + c_{1}(\lambda_1^{k} x_{1} + \hat\lambda^{k}\xi_{k1}(x, y, \varepsilon)) +  \\
\qquad + c_{2}(\lambda_2^{k} x_{2} + \hat\lambda^{k}\xi_{k2}(x, y, \varepsilon))  +
d(y - y^-)^2 + \\
\qquad + O(|y - y^-|^3 + \lambda^{k} \|x\||y - y^-| + \lambda^{2k} \|x\|^2)\;.
\\
\end{array}
\end{equation}
We shift coordinates $ x_{1new}\;=\; x_1 - x_1^+ +
\phi_k^1(\varepsilon),\; x_{2new}\;=\; x_2 - x_2^+ +
\phi_k^2(\varepsilon),\;  y_{new}\;=\; y - y^- +
\psi_k(\varepsilon)$, where $\phi_k, \psi_k = O(\lambda^{k})$,
in such a way that the right sides of (\ref{tk-1}) do not contain
constant terms for the first two equations and linear in $y_{new}$
terms for the third equation. Then (\ref{tk-1}) takes the form
\begin{equation}
\label{tk-1+}
\begin{array}{l}
\bar x_1 \; = \; a_{11}\lambda_1^{k} x_{1} + a_{12}\lambda_2^{k} x_{2} + b_1 y +
O(y^2 + \lambda^{k} \|x\| |y| + \hat\lambda^{k} \|x\|^2)  \;, \\
\\
\bar x_2  \; = \; a_{21}\lambda_1^{k} x_{1} + a_{22}\lambda_2^{k} x_{2}  + b_2 y +
O(y^2 + \lambda^{k} \|x\| |y| + \hat\lambda^{k} \|x\|^2)   \;, \\
\\
\bar y - (\hat\gamma/\gamma)^{-k}\eta_k(\bar x + x^+ + \phi_k,\bar y + y^- +
\psi_k,\varepsilon) \; = \; M_k  + \\
\qquad d \gamma^{k}y^2 + \lambda_1^{k}\gamma^{k}c_{1} x_{1} + \lambda_2^{k}\gamma^{k}c_{2}x_{2} +
\gamma^{k}O(|y|^3 + \lambda^{k} \|x\| |y|) + \hat\lambda^{k} \|x\|^2)

\end{array}
\end{equation}
where
\begin{equation}
\label{Mk1}
\begin{array}{l}
M_k\;=\;\gamma^{k}[\mu_1 + \lambda_1^{k}c_1x_1^+ + \lambda_2^{k}c_2x_2^+ + o(\lambda^k)]
\end{array}
\end{equation}
%and dots stand for coefficients tending to zero as $k\to\infty$.

Consider the third equation of (\ref{tk-1+}).
First of all, we transform its left side. Namely, we write
$\bar y -
(\hat\gamma/\gamma)^{-k}
\eta_k = \bar y +
(\hat\gamma/\gamma)^{-k}[\eta_k^0 +
\eta_k^1(\bar x,\varepsilon) + \eta_k^2(\bar y,\varepsilon) +
\eta_k^3(\bar x,\bar y,\varepsilon)]$ where
$\eta_k^1(0, \varepsilon) = 0$, $\eta_k^2(0, \varepsilon) = 0$ and
$\eta_k^3 = O(\|\bar x \bar y\|)$. Next, we transfer  constant term
$(\hat\gamma / \gamma)^{-k}\eta_k^0$ into the right side and join it to $M_k^1$;
we substitute the value of $\bar x$ due to the first two equations of (\ref{tk-1+}) into
function $\eta_k^1(\bar x, \varepsilon)$ and transfer the obtained expression
into the right side. After this, all coefficients
(in the third equation) get additions of order $O(\hat\gamma^{-k})$ and
a new linear term in $y$, $p_ky = O([\hat\gamma/\gamma]^{-k}) y$, appears.
By the shift of coordinates of the form
$(x, y) \mapsto (x,y) + O([\hat\gamma / \gamma]^{-k})$, we vanish both this linear
term and constant terms in the right sides of the first and second equations.
As the result, the left side of the third equation can be written as follows:
$$
%\begin{array}{l}
\bar y + (\hat\gamma / \gamma)^{-k}O(\bar y) +
(\hat\gamma/\gamma)^{-k}O(\|\bar x\bar y\|) =
\bar y(1 + q_k) +
(\hat\gamma/\gamma)^{-k}O(\bar y^2) +
(\hat\gamma/\gamma)^{-k}O(\|\bar x \bar y\|)
%\end{array}
$$
where $q_k = O([\hat\gamma/\gamma]^{-k})$. After this, we can write
system (\ref{tk-1+}) in the form
\begin{equation}
\label{tk-1++}
\begin{array}{l}
\bar x_1 \; = \; a_{11}\lambda_1^{k} x_{1} +  a_{12}\lambda_2^{k} x_{2} + b_1 y + O(y^2)
+ \lambda^{k} O( \|x\| |y|) + \hat\lambda^{k}O(\|x\|^2)  \;, \\
%\\
\bar x_2  \; = \; a_{21}\lambda_1^{k}x_{1} + a_{22}\lambda_2^{k} x_{2} + b_2 y +
 O(y^2) + \lambda^{k} O( \|x\| |y|) + \hat\lambda^{k}O(\|x\|^2) \;,    \\
\\
\bar y(1 + q_k) +
(\hat\gamma/\gamma)^{-k}O(|\bar y|^2) +
(\hat\gamma/\gamma)^{-k}O(\|\bar x\bar y\|)
\; = \;
M_k  + \\
\qquad + d \gamma^{k}(1 + s_k)y^2 +  c_{1} \lambda_1^{k}\gamma^{k} x_{1} + c_{2} \lambda_2^{k}\gamma^{k} x_{2} +
p_k \gamma^{k} O(\|x\|^2) + \\

\qquad + \lambda^k \gamma^{k}O(\|x\| |y|) + \gamma^{k} O(y^3)\;, \\
\end{array}
\end{equation}
where $s_k = O(\lambda^k + |\hat\gamma/\gamma|^{-k})$, $p_k = O(\hat\lambda^k + |\hat\gamma/\gamma|^{-k})$ and
new $M_k$ satisfies (\ref{Mk1}).

We perform a linear change of $x$ variables to make zero the linear in $y$ term in the second equation:
$$
\displaystyle x_{2new} = x_{2} - \frac{b_{2}}{b_1} x_1, \; x_{1new} = x_{1}\;,\; y_{new} = y.
$$

%\subsection{Proof of Theorem~\ref{th3-2}.}\label{prt3-1}
%
%Consider map  (\ref{tk-1++}) and introduce new coordinates
%$$
%\xi = x_{1} + \frac{c_{2}}{c_1} \frac{\lambda_2^{k}}{\lambda_1^{k} } x_{2} + (\lambda_{k} +
%\hat\gamma^{-k})O(\|x\|), \; x_{1new} = x_{1}\;,\; y_{new} = y,
%$$
%i.e. we take as $\xi$ the expression from the squire brackets in the third equation of
%(\ref{tk-1++}). Then (\ref{tk-1++}) is rewritten in the form

Then (\ref{tk-1++}) is rewritten in the form

\begin{equation}
\label{tk-1+-}
\begin{array}{l}
\displaystyle \bar x_1 \; = \;  b_1 y +
\lambda^{k}O(\|x\|)+ O(y^2) \;, \\
\\
\displaystyle \bar x_2 = A_{21} \lambda_1^k x_1 + A_{22}\lambda_2^k x_2 + O(y^2) + \lambda^{k} O( \|x\| |y|) + \lambda^{2k}O(\|x\|^2)\;,    \\
\\
\bar y(1 + q_k) + (\hat\gamma/\gamma)^{-k} O(|\bar y|^2 + \|\bar x\bar y\|)  \; = \; M_k +
 \lambda_1^{k}\gamma^{k}c_{1} \nu_k x_1 + \lambda_2^{k}\gamma^{k}c_{2} x_2 + d\gamma^{k}(1 + s_k)y^2 +
  \\
\qquad\qquad  p_k \gamma^{k} O(\|x\|^2) + \lambda^k \gamma^{k}O(\|x\| |y|) +

\gamma^{k}O(y^3) \;,
\end{array}
\end{equation}
where
\begin{equation}
\begin{array}{c}
\displaystyle A_{21} = \left( a_{21} - \frac{b_2}{b_1} a_{11} \right) +
\left( a_{22} - \frac{b_2}{b_1} a_{12} \right)\frac{b_2}{b_1} \frac{\lambda_2^k}{\lambda_1^k}, \\
\displaystyle A_{22} = a_{22} - \frac{b_2}{b_1} a_{12}, \\
\displaystyle \nu_k = \left( 1 + \frac{b_2c_{2}}{b_1c_1} \frac{\lambda_2^{k}}{\lambda_1^{k}} \right).
\end{array}
\label{tky}
\end{equation}

Now we will vary $\lambda_1$ and $\lambda_2$ in such a way that the value of $\nu_k$ is asymptotically small as $k \to \infty$.
This is always possible via small changes of parameter $\mu_3$ because $b_1 c_1 \neq 0$ and $b_2 c_2 \neq 0$ and $\lambda_1= -\lambda_2$ in the initial moment.
Then it is clear that
$$
\displaystyle A_{21} = \frac{J_1}{c_2 b_1} + O(\nu_k) \neq 0,
$$
where $J_1$ is given by formula (\ref{t1-3}).

Rescale the coordinates as follows
%formulae $y = \alpha_k Y \;,\; x_1 = \beta_{1k} X_1 \;,\;
%x_2 = \beta_{2k} X_2$,  where
$$
\displaystyle y  = -\frac{\gamma^{-k}(1 + q_k)}{d(1 + s_k)}\;Y \;,\;
x_1  = -\frac{b_1 \gamma^{-k}(1 + q_k)}{d(1 + s_k)}\;X_1\;,\;
x_2  = -\frac{b_1 A_{21} \gamma^{-k}\lambda_1^k (1 + q_k)}{d(1 + s_k)}\;X_2.
$$
Then system (\ref{tk-1+-}) is rewritten in the new coordinates as
follows
\begin{equation}
\label{tkres1}
\begin{array}{l}
\displaystyle \bar X_1 \; = \; Y + O(\lambda^k)\;,\; \\
\bar X_2  \; = \; X_1 + O(\gamma^{-k}\lambda^{-k})\;,\; \\
\bar Y = M_1 + M_2 X_1 + B X_2 - Y^2 + O(\gamma^{-k}\lambda^{-k})\;,
\end{array}
\end{equation}
%{tkres1}
where formulas (\ref{M1-itog}) and (\ref{M2-itog}) are valid for $M_1$, $M_2$ and $B$.

%Introduce the new coordinate
%$$
%X_{2new} = -\beta M_2 X_1 + \beta X_2,
%$$
%where $\beta = {\tilde A}^{-1}$. Then the system (\ref{tkres1}) recasts as
%\begin{equation}
%\label{tkres1-2}
%\begin{array}{l}
%\displaystyle \bar X_1 \; = \; Y + O(\lambda^k)\;,\; \\
%%
%\bar X_2  \; = X_1 + O(\gamma^{-k}\lambda_2^{-k})\;,\; \\
%%
%\bar Y
% = M_1 + \hat M_2 X_1 + J_k X_2 - Y^2 + O(\gamma^{-k}\lambda_2^{-k})\;,
%%
%\end{array}
%\end{equation}
%where  $J_k = J_1\lambda_1^{k}\lambda_2^k\gamma^{k}(1 + \dots) $.  We obtain also that
%$$
%\displaystyle M_2 = \left( b_1c_1 + b_2c_{2} \frac{\lambda_2^{k}}{\lambda_1^{k}}
%\right)\lambda_1^k\gamma^k(1+\dots)
%$$

It is obvious that system (\ref{tkres1}) is asymptotically close to (\ref{tk-itog}) when $k \to \infty$. $\;\;\;\Box$.

\end{document}